\documentclass[a4paper]{article}
\usepackage[latin1]{inputenc}
\usepackage{amscd}
\usepackage[english]{babel}
\usepackage{amsmath,amssymb,amsfonts}
\usepackage{syntonly}

\usepackage{color}
\usepackage[outerbars,color]{changebar}
\usepackage{ifthen}

\newboolean{usecolor}
\setboolean{usecolor}{true}

\ifthenelse{\boolean{usecolor}}
{
  \definecolor{colore}{cmyk}{0,1,0.6,0}
  \definecolor{coloregen}{cmyk}{0.7,0,1,0}
  \definecolor{coloresimo}{cmyk}{1,0.6,0,0}
}
{
  \definecolor{colore}{cmyk}{0,0,0,1}
  \definecolor{coloregen}{cmyk}{0,0,0,1}
  \definecolor{coloresimo}{cmyk}{0,0,0,1}
}

\newcommand{\notagen}[1]
{\-\marginpar[\raggedleft\tiny \textcolor{coloregen}{#1 G}]{\raggedright\tiny \textcolor{coloregen}{#1 G}}}

\newcommand{\notasimo}[1]
{\-\marginpar[\raggedleft\tiny \textcolor{coloresimo}{#1 S}]{\raggedright\tiny \textcolor{coloresimo}{#1 S}}}

\newcommand{\rx}{\mathcal{R}_X}

\newcommand{\Z}[0]{\mathbb{Z}}

\newcommand{\C}[0]{\mathbb{C}}

\newcommand{\Cal}[1]{\mathcal{#1}}

\newcommand{\eq}[1][r]
       {\ar@<-3pt>@{-}[#1]
        \ar@<-1pt>@{}[#1]|<{}="gauche"
        \ar@<+0pt>@{}[#1]|-{}="milieu"
        \ar@<+1pt>@{}[#1]|>{}="droite"
        \ar@/^2pt/@{-}"gauche";"milieu"
        \ar@/_2pt/@{-}"milieu";"droite"}
\newcommand{\imm}[1][r] {\ar@{^{(}->}[#1]}

\newtheorem{df}{Definition}[section]
\newtheorem{teo}{Theorem}

\newtheorem{lem}[df]{Lemma}
\newtheorem{conj}[df]{Conjecture}

\newtheorem{rmk}[df]{Remark}

\newcommand{\wdt  }{\widetilde}

\newcommand{\bW}{W}

\newcommand{\lgr}{\longrightarrow}

\title{The integer cohomology of toric Weyl arrangements}
\author{Simona {\sc Settepanella}\footnote{LEM, Scuola Superiore Sant'Anna, Pisa, Italy. s.settepanella@sssup.it ( Thanks to financial support from the European Commission 6th FP (Contract CIT3-CT-2005-513396), Project: DIME - Dynamics of Institutions and Markets in Europe) }}
\date {}

\begin{document}

\maketitle
\begin{abstract}
A toric arrangement is a finite set of hypersurfaces in a complex torus, every hypersurface being the kernel of a character. In the present paper we prove that if $\Cal T_{\wdt W}$ is the toric arrangement defined by the \textit{cocharacters} lattice of a Weyl group $\wdt W$, then the  integer cohomology of its complement is torsion free.
\end{abstract}

\begin{center}
\end{center}

\begin{center}
\end{center}

\section*{Introduction}

Let $T=(\C^*)^n$ be a complex torus and $X \subset Hom(T,\C^*) $ be a finite set of
characters of $T$. The kernel of every $\chi \in X$ is a hypersurface
of $T$:
\begin{equation*}
H_{\chi}:=\{t \in T \, \mid \, \chi(t)=1\}.
\end{equation*}
Then $X$ defines on T the \textit{toric arrangement}:
\begin{equation*}
\mathcal{T}_X:=\{ H_{\chi} , \chi \in X\}.
\end{equation*}
Let $\rx$ be the \textit{complement } of the arrangement:
\begin{equation*}
\rx:= T \setminus \bigcup_{\chi \in X} H_{\chi}
\end{equation*}

The geometry and topology of $\rx$ have been studied by many authors, 
see for instance \cite{L1}, \cite{L2}, \cite{li}, \cite{ERS}, \cite{Mo2} and \cite{Mo3}. 
In particular Looijenga (see \cite{Lo}) and De Concini and Procesi (see \cite{DP}) computed the De Rham cohomology of $\rx$ and, recently, Moci and Settepanella (see \cite{MoSe}) described a regular CW-complex homotopy equivalent to $\rx$. This complex is similar to the one introduced by Salvetti (see \cite{Sa1}) for the complement of hyperplane arrangements. 

If $\Cal T_{\wdt W}$ is the toric arrangement associated to an affine Weyl group $\wdt W$, the complex $T(\wdt W)$ homotopic to the complement
$$
\Cal R_{W}:= T \setminus \bigcup_{H \in T_{\wdt W}} H
$$
admits a very nice description which generalizes a construction introduced in \cite{Sa} and \cite{boss3}. 
In their paper Moci and Settepanella conjectured that the integer cohomology of $T(\wdt W)$ 
(equivalently $\Cal R_W$) is torsion free. Hence it coincides with the De Rham cohomology described in \cite{DP} and
it is known since the Betti numbers can be easily computed using results in \cite{Mo1}.

\bigskip

In the present paper we prove this conjecture generalizing to toric arrangements a well known result  for hyperplane ones. Indeed Arnol'd proved that the integer cohomology of braid arrangement is torsion free in 1969 (see \cite{Arno}).

\bigskip

In order to prove it we use a filtration introduced in \cite{DPS} and generalized to braid arrangements in \cite{simo1} (see subsection \ref{seziofiltra}). 

In Section \ref{sezione2} we prove that the above filtration involves complexes with torsion free cohomology. While in Section \ref{sezione3} we rewrite it for toric arrangements and we prove the main result of the paper:

\begin{teo}The integer (co)-homology of the complement $\Cal R_W$  is torsion free.
\end{teo}

\bigskip

\paragraph{Acknowledgements}
I wish to thank Filippo Callegaro and Luca Moci for very useful conversations we had while I was revising the present paper.

\section{Notations and recalls}\label{S-complex} 

In this section we recall basic construction about affine and toric arrangements coming from 
Coxeter systems.

\subsection{Salvetti's complex for Coxeter arrangements}

Let $(W,S)$ be the Coxeter system associated to the finite reflection group $W$ and
$$\mathcal{A}_{W}=\{H_{w s_i w^{-1}} \mid w
\in W \mbox{ and } s_i \in S\}$$
the arrangement in $\C^n$ obtained by complexifying the reflection hyperplanes of
$W$, where, in a standard way, the hyperplane
$H_{w s_i w^{-1}}$ is simply the hyperplane fixed by the reflection
$w s_i w^{-1}$.

\bigskip

It is well known (see, for instance, \cite{boss3} \cite{Sa} ) that the $k$-cells of Salvetti's complex
$C(W)$ for arrangements $\Cal{A}_{W}$ are of the form
$E(w,\Gamma)$ with
$\Gamma \subset S$ of cardinality $k$ and $w \in  W$. 

While the integer boundary map can be expressed as follows:

\begin{equation}
\begin{split}
\partial_k(E(w,\Gamma)) =\\
&\sum_{s_j \in \Gamma}
\sum_{\beta\in W^{\Gamma\setminus\{ s_j \}}_{\Gamma}}(-1)^{l(\beta)+\mu
(\Gamma,s_j)} E(w\beta,\Gamma\setminus\{ s_j\})
\end{split}
\end{equation}
where $W_{\Gamma}$ is the group generated by $\Gamma$, 
$$W^{\Gamma\setminus\{\sigma\}}_{\Gamma}=\{w \in W_{\Gamma} : l(ws) > l(w) \forall s \in \Gamma \setminus\{\sigma\} \}$$ 
and $\mu(\Gamma, s_j)=\sharp\{s_i \in \Gamma | i \leq j \}$.
Here $l(w)$ stands for the length of $w$.

\begin{rmk}Instead of the co-boundary operator we prefer to describe its dual,
i.e. we define the boundary of a $k$-cell $E(w,\Gamma )$ as a linear combination of
the $(k-1)$-cells which have $E(w,\Gamma )$ in theirs co-boundary, with the
same coefficient of the co-boundary operator.
We make this choice since the boundary operator has a nicer description than
co-boundary operator in terms of the elements of $W$.
\end{rmk}

This description holds also for Coxeter systems $(\wdt W,\wdt S)$ associated to 
Weyl groups $\wdt W$.

\subsection{A filtration for the complex $(C(W),\partial)$}\label{seziofiltra}

It's known (see \cite{Bou}) that for all $\Gamma \subset S$ the group $W$ splits as
\begin{equation*}
W=W^{\Gamma}W_{\Gamma}
\end{equation*}
with
\begin{equation}\label{wgamma}
W^{\Gamma}=\{w^{\Gamma} \in W \mid l(w^{\Gamma}s_i)>l(w^{\Gamma}) \mbox{ for all } s_i \in W_{\Gamma}\}.
\end{equation}
If $w=w^{\Gamma}w_{\Gamma}$, $w^{\Gamma} \in W^{\Gamma}$ and $w_{\Gamma} \in W_{\Gamma}$, then $l(w\beta)=l(w^{\Gamma})+l(w_{\Gamma}\beta)$ $\forall \beta \in \bW_{\Gamma}$ and the boundary map verifies
$$
\partial(E(w,\Gamma)) = 
w^{\Gamma}.\partial(E(w_{\Gamma},\Gamma)). 
$$

In  \cite{simo1} (see also \cite{DPS})  author defines a map of complexes
\begin{equation*}
i_m := i : \bigoplus_{j=1}^{m_1} C(W_{S_{m-1}}) \lgr C(\bW)
\end{equation*}
 as follows 
\begin{equation*}
\begin{split}
&i(j.E(w_{S_{m-1}},\Gamma))=i(\bW^{S_{m-1}}(j).E(w_{S_{m-1}},\Gamma))=\\
&i(w^{S_{m-1}}.E(w_{S_{m-1}},\Gamma)=w^{S_{m-1}}.i(E(w_{S_{m-1}},\Gamma))=
w^{S_{m-1}}.E(w_{S_{m-1}},\Gamma)=E(w,\Gamma).
\end{split} 
\end{equation*}
Where $m_1$ is the cardinality of $W^{S_{m-1}}$, 
$w^{S_{m-1}}=\bW^{S_{m-1}}(j)$ its $j$-th element in a fixed order 
and $S_{h}=\{s_1,\cdots,s_h\} \subset S=\{s_1,\cdots ,s_m\}$.

The cokernel of the map $i$ is the complex $F^1_m(\bW)$
having as basis all 
$E(w,\Gamma_1)$ for $w \in \bW$ and $\Gamma_1 \subset S$ s.t. $s_m \in \Gamma_1$.

She iterates this construction getting maps
\begin{equation*}
\begin{split}
&i_m[k] := i : \bigoplus_{j=1}^{m_1\cdots m_{k+1}}C(\bW_{S_{m-k-1}})[k] 
\lgr F_m^k(\bW),\\
&i(w^{S_{m-k-1}}.(E(w_{S_{m-k-1}},\Gamma)))=
w^{S_{m-k-1}}.i(E(w_{S_{m-k-1}},\Gamma)) \\
&=E(w,\Gamma \cup \{s_m,\cdots,s_{m-k+1}\})
\end{split} 
\end{equation*}  

Each $i_m[k]$ gives rise to the exact sequence of complexes
\begin{equation} \label{eqn:suc.es.}
0 \lgr \bigoplus_{j=1}^{m_1\cdots m_{k+1}}C(\bW_{S_{m-k-1}})[k] \stackrel{i}{\lgr} 
F^k_m(\bW) \stackrel{j}{\lgr} F^{k+1}_{m}(\bW) \lgr 0.
\end{equation}

It is possible to filter the complex $F^0_m(W)=C(W)$ in a similar way through maps:
\begin{equation}\label{eqn4}
\begin{split}
&i^m[k] := i : \bigoplus_{j=1}^{m^1\cdots m^{k+1}}C(\bW_{S^{k+1}})[k] 
\lgr F^k_m(\bW),\\
&i(w^{S^{k+1}}.(E(w_{S^{k+1}},\Gamma)))=
w^{S^{k+1}}.i(E(w_{S^{k+1}},\Gamma)) =\\
&=E(w,\{s_1,\cdots,s_{k} \} \cup \Gamma\})
\end{split} 
\end{equation}  
for $0 \leq k \leq m$, $S^k=\{s_{k+1},\cdots,s_m\}$ and $m^i$ the cardinality of $W_{S^{i-1}}^{S^i}$.

\subsection{Salvetti's complex for toric Weyl arrangements}\label{seztor}

Let $\Phi$ be a root system, $\langle\Phi^\vee\rangle$ be the
lattice spanned by the coroots, and $\Lambda$ be its dual
lattice (which is called the \emph{cocharacters} lattice).
Then we define a torus $T=T_\Lambda$ having $\Lambda$
as group of characters.

If $\widetilde{W}$ is the affine Weyl group associated to $\Phi$,
we can regard $\Lambda$ as a subgroup of $\widetilde{W}$, acting by translations. It is well known that $\widetilde{W}/\Lambda \simeq W$, where $W$ is the finite reflection group associated to $\widetilde{W}$. As a consequence, the toric Weyl arrangement
can be described as:
$$
\Cal T_{\widetilde{W}}=\{H_{[w] s_i [w^{-1}]} \mid w \in W \mbox{ and } s_i \in \wdt S \}
$$
where two hypersurfaces $H_{[w] s_i [w^{-1}]}$ and $H_{[\overline{w}] s_i [\overline{w}^{-1}]}$ are equal if and only if
there is a translation $t \in \Lambda$ such that $tw s_i (tw)^{-1}=\overline{w} s_i \overline{w}^{-1}$, i.e. $\overline{w}=tw$.

\bigskip

In \cite{MoSe} authors prove that the complement 
$$
\Cal R_W := T \setminus \bigcup_{H \in T_{\wdt W}} H
$$
has the same homotopy type of a CW-complex $T(\wdt W)$ which admits a description similar to $C(W)$.

Indeed the $k$-cells of $T(\wdt W)$ correspond to elements $E([w],\Gamma)$ where $[w] \in \widetilde{W} / \Lambda \simeq W$ is an equivalence class  with one and only one representative $w \in W$ and
$\Gamma=\{s_{i_1},\ldots , s_{i_k}\}$ is a subset of cardinality $k$ in $\wdt S$.

The integer boundary operator is 
\begin{equation}\label{bordo}
\begin{split}
\partial_k(E([w],\Gamma)) =\\
&\sum_{\sigma\in \Gamma}
\sum_{\beta\in W^{\Gamma\setminus\{\sigma\}}_{\Gamma}}(-1)^{l(\beta)+\mu
(\Gamma,\sigma)} E([w\beta],\Gamma\setminus\{ \sigma \}).
\end{split}
\end{equation}

Let $\Gamma \subset \wdt S$ be a proper subset and $W_{\Gamma}$ be the finite reflection group generated by $\Gamma$. The group
$$
(\wdt W / \Lambda)_{\Gamma}=\{[w] \in \wdt W / \Lambda \mid w \in W_{\Gamma}\} \simeq W_{\Gamma}
$$
is a well defined subgroup of $\wdt W / \Lambda$. As in the finite case, we get 
$$\wdt W / \Lambda = (\wdt W / \Lambda)^{\Gamma}(\wdt W / \Lambda)_{\Gamma}$$
and the toric boundary map verifies
$$
\partial(E([w],\Gamma)) = 
[w^{\Gamma}].\partial(E([w_{\Gamma}],\Gamma)) 
$$
where $[w^{\Gamma}] \in (\wdt W / \Lambda)^{\Gamma}$,  $[w_{\Gamma}] \in 
(\wdt W / \Lambda)_{\Gamma}$ and 
$[w]=[w^{\Gamma}] [w_{\Gamma}]=[w^{\Gamma}w_{\Gamma}].$

\bigskip

Let us remark that $(\wdt W / \Lambda)_{\Gamma}$ is isomorphic to a subgroup 
of $W$ which is not, in general, a parabolic one. In these cases the set $(\wdt W / \Lambda)^{\Gamma}$ 
doesn't admit a description similar to the one in (\ref{wgamma}). 

Our main interest in the sequel of this paper is to construct a filtration for $T(\wdt W)$ similar to the one in subsection \ref{seziofiltra}. 
Also if it is not necessary to know an explicit description of $(\wdt W / \Lambda)^{\Gamma}$ in order to filter the complex $T(\wdt W)$, nevertheless we believe that it would be useful to know a little bit more about it to have a better understanding of our construction.
In particular, if $\wdt S=\{s_0,\ldots ,s_m\}$, we are interested in the cases in which $\Gamma = \{s_k,\ldots ,s_m\}$ or $\Gamma=\{s_0,\ldots s_h\}$.

It is a simple remark that, if $s_0 \notin \Gamma$, then 
$(\wdt W / \Lambda)_{\Gamma} \simeq W_{\Gamma}$ is a parabolic subgroup of $W$. 
While the case $s_m \notin \Gamma$ is a little bit more complicated. 
Since to remove $s_0$ or $s_m$ is perfectly symmetric for $\wdt W= \wdt A_m, \wdt C_m, \wdt D_m, \wdt E_6, \wdt E_7$, then in these cases we always get that $(\wdt W / \Lambda)_{\Gamma} \simeq W_{\Gamma}$ is a parabolic subgroup of $W$.  Hence in the above situations 
$(\wdt W / \Lambda)^{\Gamma} \simeq W^{\Gamma}$ admits a description as in (\ref{wgamma}). 

Otherwise $W_{\wdt S\setminus \{s_m\}}$ is still a finite reflection group but it is not of type $W$. For example if $\wdt W=\wdt B_m$ then  $W_{\wdt S\setminus \{s_m\}}= D_{m}$ which is not 
$B_m$.  In these cases if $\Gamma \subset \wdt S$ 
is a given subset with $s_m \notin \Gamma$ and $s_0 \in \Gamma$,
 then $(\wdt W / \Lambda)_{\Gamma} \simeq W_{\Gamma}$ is a parabolic subgroup of $W_{\wdt S\setminus \{s_m\}}$ and, by \cite{Mo1}, we have exactly
$$
\frac{\mid W \mid}{\mid W_{\wdt S \setminus \{s_m \} } \mid}
$$
copies of $W_{\wdt S\setminus \{s_m\}}$ in $W$. 

Let $W^{\prime}$ be the subgroup of $W$ such that
$W^{\prime} \simeq W_{\wdt S\setminus \{s_m\}} \simeq (\wdt W / \Lambda)_{\wdt S\setminus \{s_m\}}$ then $W^{\wdt S\setminus \{s_m\}}$ will denote the subset of $W$ such that 
$W=W^{\wdt S \setminus \{s_m\}}W^{\prime}$ and  we get
$$
(\wdt W / \Lambda)^{\Gamma} \simeq W^{\wdt S \setminus \{s_m \}}W_{\wdt S \setminus \{s_m \}}^{\Gamma}
$$
where $W_{S \setminus \{s_m \}}^{\Gamma}$ is the subset of $W_{\wdt S \setminus \{s_m \}}$ described in (\ref{wgamma}).

\section{The cohomology of complexes $F^k_n(\bW)$}\label{sezione2}

It is well known that the integer homology, and hence cohomology, of complexes $C(\bW)$ is 
torsion free, while the (co)-homology $H^*(F^k_n,\Z)$ is not known. 
In this section we will prove that it is torsion free.

\bigskip

As above we will consider the boundary map instead of the (co)-boundary one. \\
The exact sequences (\ref{eqn:suc.es.}) give rise to the corresponding long 
exact sequences in homology
\begin{equation}\label{succomologia} 
\begin{split}
\cdots \lgr H_{*+1}(F^{k}_m(\bW),\Z) 
&\stackrel{\Delta_*}{\lgr} \bigoplus_{j=1}^{m_1\cdots m_{k}}H_{*-k}(C(\bW_{S_{m-k}}),\Z) 
\stackrel{i_*}{\lgr}\\
&\stackrel{i_*}{\lgr} H_{*}(F^{k-1}_m(\bW),\Z) \stackrel{j_*}{\lgr} H_{*}(F^{k}_m(\bW),\Z) \lgr 
\cdots 
\end{split} 
\end{equation}
where the map $\Delta_*$ is induced by the map on complexes:
\begin{equation} \label{delta}
\begin{split}
&\Delta: F^k_m(\bW) \lgr \bigoplus_{j=1}^{m_1\cdots m_{k}}C(\bW_{S_{m-k}})\\
&\Delta(E(w,\Gamma \cup S^{m-k}))=
\sum_{\beta \in W^{\Gamma \cup S^{m-k+1}}_{\Gamma \cup S^{m-k}}}(-1)^{l(\beta)}E(w\beta, \Gamma).\\
\end{split}
\end{equation}

%
%
To simplify notation from now on we will use 
$$l=m-k-1$$ 
and $\bigoplus C(\bW_{S_{m-k}})$ instead of
$\bigoplus_{j=1}^{m_1\cdots m_{k}}C(\bW_{S_{m-k}})$ since the number of copy 
$\prod_{i=1}^k m_i$ is completely determined by $S_{m-k}$. 

We have the following theorem.

\begin{teo}\label{teofmk}The integer (co)-homology of complexes $F^k_m(\bW)$ is 
torsion free for all $k \leq m$.
\end{teo}

We need the following key Lemma.

\begin{lem}\label{chiave}Let $v \in F^k_m(\bW)$ be a boundary then one of the following occurs:

\bigskip

i) $v \in i(\bigoplus C(\bW_{S_{l}})[k])$ 

ii) $v \in F^k_m(\bW) \setminus i(\bigoplus C(\bW_{S_{l}})[k])$ 

iii) $v$ is a sum of two boundaries $v^{\prime} \in i(\bigoplus C(\bW_{S_{l}})[k])$ and 
$v^{\prime \prime} \in F^k_m(\bW) \setminus i(\bigoplus C(\bW_{S_{l}})[k])$.
\end{lem}

\textbf{Proof.} By construction any chain $v \in F^k_m(\bW)$ is a sum of two chains 
$$v=v^{\prime} + v^{\prime \prime}$$
the first one in $i(\bigoplus C(\bW_{S_{l}})[k])$ and the second one in
$F^k_m(\bW) \setminus i(\bigoplus C(\bW_{S_{l}})[k])$. 
Let $v$ be a boundary. If $v^{\prime}$ ($v^{\prime \prime}$) is zero then $ii)$ ($i) )$ follows.

Let $v^{\prime}$ and $v^{\prime \prime}$ both not zero.
Ordering in a suitable way rows and columns of the boundary matrix, we get a block matrix as follows:
\begin{displaymath}
\left[
\begin{array}{cc}
\bigoplus  i(\partial C(\bW_{S_{l}})[k]) & B_1    \\
          0             & B_2 
\end{array}
\right] \quad \mbox{and} \quad
\left[
\begin{array}{c}
B_1 \\
B_2
\end{array}
\right]
=\partial (F^{k}_m(\bW) \setminus i(\bigoplus C(\bW_{S_{l}})[k]).
\end{displaymath}
Then we can diagonalize the matrix by row and column operations in such a way that the rows of 
the first (second) block are combined only with rows in the same block.\\
As consequence any element $v$ which is in the boundary is written as a sum of two boundaries , one obtained by combinations of row in the first block, i.e. a combination of elements in 
$i(\bigoplus C(\bW_{S_{l}})[k])$, and the second one by elements in $F^k_m(\bW) \setminus i(\bigoplus C(\bW_{S_{l}})[k])$. $\qquad$ 
$\square$

\bigskip

\begin{rmk}If $v^{\prime}$ and $v^{\prime \prime}$ are boundaries in $F^k_m(\bW)$
as in the above Lemma, then 
$v^{\prime} \in i( \bigoplus \partial C(\bW_{S_{l}})[k])$ while, obviously, $v^{\prime \prime}$ is a linear combination of elements in $F^k_m(\bW) \setminus i(\bigoplus C(\bW_{S_{l}})[k])$, but it is not in its boundary. 
\end{rmk}

\bigskip

\textbf{Proof of Theorem \ref{teofmk}}  The integer cohomology of the complex 
$F^0_m(\bW) = C(\bW)$ is torsion free.
By induction let us assume that $H^*(F^{k-1}_m (\bW), \Z)$, and hence 
$H_*(F^{k-1}_m(\bW), \Z)$, are torsion free.

\bigskip

As the sequence (\ref{succomologia}) is exact and $H_*(C(\bW_{S_{m-k}}),\Z)$ and $H_*(F^{k-1}_m (\bW), \Z)$ are torsion free, then
$H_*(F^k_m (\bW), \Z)$ (and hence $H^*(F^k_m (\bW), \Z)$) is torsion free if and only if the 
image of $i_*$ doesn't give rise to $p$-torsion for $p \in \Z$, i.e. 
\begin{equation*}
p[v] \in i_*(\bigoplus_{j=1}^{m_1\cdots m_{k}}H_*(C(\bW_{S_{m-k}}),\Z) ) \Longleftrightarrow
[v] \in i_*(\bigoplus_{j=1}^{m_1\cdots m_{k}}H_*(C(\bW_{S_{m-k}}),\Z) ).
\end{equation*}
Let $[v] $ be a generator in the free module $H_*(F^{k-1}_m(W), \Z)$.
By construction 
$$[v]=z^{\prime} + z^{\prime \prime} + \partial_*(F^{k-1}_m(\bW))$$
for $z^{\prime} \in i(\bigoplus C(\bW_{S_{l}})[k])$ and $z^{\prime \prime} \in F^{k-1}_m(\bW) \setminus i(\bigoplus C(\bW_{S_{l}})[k])$.\\
Let us assume $$p[v]=pz^{\prime} +pz^{\prime \prime} + \partial_*(F^{k-1}_m(\bW)) \in i_*(\bigoplus H_*(C(\bW_{S_{m-k}}),\Z) ).$$ 
Then $p[v]$ has at list one representative in  the image $i(\bigoplus C(\bW_{S_{l}})[k])$ and hence there is an element
$$\omega= \omega^{\prime} + \omega^{\prime \prime} \in  \partial_*(F^{k-1}_m(\bW))$$
such that $pz^{\prime}+pz^{\prime \prime}+ \omega \in i(\bigoplus C(\bW_{S_{l}})[k])$, i.e.
$\omega^{\prime} \in  i(\bigoplus C(\bW_{S_{l}})[k])$ and 
$\omega^{\prime \prime}= -pz^{\prime \prime}$.

\bigskip 

 By Lemma \ref{chiave} we get that $-\omega^{\prime \prime} = pz^{\prime \prime} \in \partial_*(F^{k-1}_m(\bW))$ and hence
 $z^{\prime \prime} \in \partial_*(F^{k-1}_m(\bW))$ since $H_*(F^{k-1}_m(\bW))$ has no torsion by inductive hypothesis. 
Then $$[v]=z^{\prime} + z^{\prime \prime} + \partial_*(F^{k-1}_m(\bW))=  z^{\prime} + \partial_*(F^{k-1}_m(\bW))$$
i.e.  $[v] \in i_*(\bigoplus H_*(C(\bW_{S_{m-k}}),\Z) )$ $\qquad$ $\square$

\bigskip

\begin{rmk} Obviously Theorem \ref{teofmk} holds also for complexes $F^k_m(\bW)$ obtained 
filtering with the inclusions in (\ref{eqn4})
\end{rmk}

An important consequence of the above theorem is that maps $\Delta_*$ are map between finitely generated free modules such that 
\begin{equation*}
p[v] \in \Delta_*(H_*(F^k_m(W),\Z) ) \Longleftrightarrow
[v] \in \Delta_*(H_*(F^k_m(W),\Z) ).
\end{equation*}
A map between two free modules which satisfies the above condition will be called
a \textit{one-free} map and it can be diagonalized as:
\begin{displaymath}
\left[
\begin{array}{cc}
I & 0 \\
0  & 0
\end{array}
\right] 
\end{displaymath}
where $I$ is the identity matrix. It is a simple remark that composition of one-free maps is still a one-free map.
This will be useful in the next section.

\section{The integer cohomology of $\Cal R_W$}\label{sezione3}

In this section we prove that the (co)-homology of $T(\wdt W)$ (i.e. $\Cal R_W$) 
is torsion free. In order to do it we construct a filtration of $T(\wdt W)$ similar to the one of $C(W)$. 

\subsection{A filtration for the complex $(T(\wdt \bW),\partial)$}
Let $\wdt S=\{s_{0},\cdots,s_m\}$ be the system of generators of $\wdt W$ and $W$ the finite group associated. We will keep the notation $S^k=\{s_{k+1}, \ldots ,s_m\} \subset \wdt S$ while we introduce the new one $\wdt S_h =\{s_0,\ldots ,s_h\} \subset \wdt S$. 

Let us consider  the natural inclusion
\begin{equation*}
i_m := i : \bigoplus_{j=1}^{m_1} C(\bW_{\wdt S_{m-1}}) \lgr  T(\wdt \bW), 
\end{equation*}
defined as:
\begin{equation*}
\begin{split}
&i(j.E(w_{\wdt S_{m-1}},\Gamma))=i(W^{\wdt S_{m-1}}(j).E(w_{\wdt S_{m-1}},\Gamma))=\\
&i(w^{\wdt S_{m-1}}.E(w_{\wdt S_{m-1}},\Gamma))= [w^{\wdt S_{m-1}}].E([w_{\wdt S_{m-1}}],\Gamma)=E([w],\Gamma)
\end{split} 
\end{equation*}
where $m_1$ is the cardinality of the set $W^{\wdt S_{m-1}}=W^{\wdt S \setminus \{s_m\}}$ 
defined in subsection \ref{seztor} and $w^{\wdt S_{m-1}}=W^{S_{m-1}}(j)$ its $j$-th element in a fixed order. 

Let us remark that $m_1$ could be also equal to $1$ depending on the 
type of $\wdt W$ as seen in subsection \ref{seztor}.

The cokernel of the map $i$ is the toric complex $F^1_m(\wdt \bW)$ having as basis all 
$E([w],\Gamma_1)$ for $w \in \bW$ and $\Gamma_1 \subset \wdt S$ with $\mid \Gamma_1 \mid \leq m$ s.t. $s_m \in \Gamma_1$.

We can iterate this construction getting maps
\begin{equation*}
\begin{split}
&i_m[k] := i : \bigoplus_{j=1}^{m_1\cdots m_{k+1}}C(\bW_{\wdt S_{l}})[k] 
\lgr F^k_m(\wdt \bW),\\
&i(w^{\wdt S_{l}}.E(w_{\wdt S_{l}},\Gamma))=
[w^{\wdt S_{l}}].E([w_{\wdt S_{l}}],\Gamma) =E([w], \Gamma \cup S^{m-k})
\end{split} 
\end{equation*}  
with $l=m-k-1$.

Each $i_m[k]$ gives rise to the exact sequence of complexes
\begin{equation} \label{eqn:suc.es.tor}
0 \lgr \bigoplus_{j=1}^{m_1\cdots m_{k+1}}C(\bW_{\wdt S_{l}})[k] \stackrel{i}\lgr 
F^k_m(\wdt \bW) \stackrel{j}\lgr F^{k+1}_{m}(\wdt \bW) \lgr 0.
\end{equation}

In a similar way we can filter using  the inclusion:
\begin{equation*}
\begin{split}
&i^m := i : C(\bW_{S^0}) \lgr  T(\wdt \bW), \\
&i(E(w,\Gamma))= E([w],\Gamma).
\end{split} 
\end{equation*}

Here $C(\bW_{S^0})$ is the classical Salvetti's complex for the finite reflection group 
$W_{S^0}=W$. The cokernel of the map $i$ is the toric complex $F^1_m(\wdt \bW)$ having as basis all 
$E([w],\Gamma_1)$ for $w \in \bW$ and $\Gamma_1 \subset \wdt S$ with $\mid \Gamma_1 \mid \leq m$ s.t. $s_0 \in \Gamma_1$.

We can iterate this construction getting maps
\begin{equation*}
\begin{split}
&i^m[k] := i : \bigoplus_{j=1}^{m^1\cdots m^{k}}C(\bW_{S^{k}})[k] 
\lgr  F^k_m(\wdt \bW),\\
&i(w^{S^{k}}.(E(w_{S^{k}},\Gamma)))=
[w^{S^{k}}].i(E([w_{S^{k}}],\Gamma)) =E([w],\Gamma \cup \wdt S_{k-1}).
\end{split} 
\end{equation*}  

Each $i^m[k]$ gives rise to the exact sequence of complexes
\begin{equation} 
0 \lgr \bigoplus_{j=1}^{m^1\cdots m^{k}}C(\bW_{S^{k}})[k] \stackrel{i}\lgr 
F^k_m(\wdt \bW) \stackrel{j}\lgr F^{k+1}_{m}(\wdt \bW) \lgr 0.
\end{equation}

\subsection{Computation of integer cohomology}

The exact sequences (\ref{eqn:suc.es.tor}) give rise to the corresponding long 
exact sequences in homology
\begin{equation*}
\begin{split}
\cdots \lgr H_{*+1}( F^{k+1}_m(\wdt \bW),\Z) 
&\stackrel{\wdt \Delta_*}\lgr \bigoplus_{j=1}^{m_1\cdots m_{k+1}}H_{*-k}(C(\bW_{\wdt S_{l}}),\Z) 
\stackrel{i_*}\lgr\\
&\stackrel{i_*}\lgr H_{*}(F^{k}_m(\wdt \bW),\Z) \stackrel{j_*}
\lgr H_{*}(F^{k+1}_m(\wdt \bW),\Z) \stackrel{\wdt \Delta_*}\lgr 
\cdots .
\end{split} 
\end{equation*}

The map $\wdt \Delta_*$ is the one induced by maps on complexes:
\begin{equation} \label{deltatilda}
\begin{split}
&\wdt \Delta: F^{k+1}_m(\wdt \bW) \lgr \bigoplus_{j=1}^{m_1\cdots m_{k+1}}C(\bW_{\wdt S_{l}})\\
&\wdt \Delta(E([w], \Gamma \cup S^{l} ))=
\sum_{\beta \in W^{\Gamma \cup S^{l+1}}_{\Gamma \cup S^{l}}}(-1)^{l(\beta)}E([w\beta], \Gamma).\\
\end{split}
\end{equation}

If  $H_{*}( F^{k+1}_m(\wdt \bW),\Z)$ are torsion free, then $H_{*}( F^{k}_m(\wdt \bW),\Z)$ 
are torsion free if and only if the maps $\wdt \Delta_*$ are one-free maps, i.e. if a generator $[u] \in \bigoplus_{j=1}^{m_1\cdots m_{k+1}}H_{*-k}(C(\bW_{\wdt S_{l}}),\Z)$ is such that
$p[u] \in Im \wdt \Delta_*$ for an integer $p \in \Z$, then $[u]  \in Im \wdt \Delta_*$.
We will prove it through an inductively argument.

\bigskip

When $k=m-1$ we get the last long exact sequence in homology
\begin{equation*} 
\begin{split}
0 \lgr   \bigoplus_{j=1}^{m_1\cdots m_{m}}&H_{1}(C(\bW_{\wdt S_0}),\Z) \stackrel{i_*}\lgr 
H_{m}(F^{m-1}_m(\wdt \bW),\Z) \stackrel{j_*}\lgr H_{m}(F^{m}_m(\wdt \bW),\Z) 
\stackrel{\wdt \Delta_*}\lgr \\
& \stackrel{\wdt \Delta_*}\lgr \bigoplus_{j=1}^{m_1\cdots m_{m}}H_{0}(C(\bW_{\wdt S_0}),\Z) 
\stackrel{i_*}\lgr H_{m-1}(F^{m-1}_m(\wdt \bW),\Z) \lgr 0 .
\end{split} 
\end{equation*}

As in the affine case, we drop the indices $m_i$ from the sum $\bigoplus$ when no misunderstanding is possible.

\bigskip

The integer homology for affine arrangements is torsion free and
$$H_{m}(F^{m}_m(\wdt \bW),\Z)=F^m_m(\wdt \bW) \simeq F^m_m(\bW) = H_{m}(F^m_m(\bW),\Z)$$ are the free modules generated by 
$E([w],S^{0})=E([w],S) \simeq E(w,S)$. Moreover, by definition,
the map $\wdt \Delta$ acts on $F^m_m(\wdt \bW)$ as $\Delta$ on $F^m_m(\bW)$.

Hence, if $C(W_{\emptyset})$ denotes the complex generated by the $0$-cell $E(1,\emptyset)$, we get the following commutative diagram in homology:

\begin{equation} \begin{array}{ccc}
 H_{m}(F^m_m(\bW),\Z)&\xrightarrow{\Delta_*} &  \bigoplus_{j=1}^{\sharp W}H_{0}(C(\bW_{\emptyset}),\Z) \\
\wr \mid & & \downarrow\scriptstyle{i_*}  \\
H_{m}(F^{m}_m(\wdt \bW),\Z)& \xrightarrow{\wdt \Delta_*} & \bigoplus_{j=1}^{m_1\cdots m_{m}}H_{0}(C(\bW_{\wdt S_0}),\Z) 
\end{array}
\end{equation}
induced by the corresponding maps on complexes.
Then, if $k=m-1$, $\wdt \Delta_*$ is one-free as composition of two one-free maps $\Delta_*$ and $i_*$ and 
$H_{m-1}(F^{m-1}_m(\wdt \bW),\Z)$ is torsion free. This provide the base of induction.
 
We remark that $H_{m}(F^{m-1}_m(\wdt \bW),\Z)$ is torsion free since the map 
$$0 \lgr \bigoplus_{j=1}^{m_1\cdots m_{m}}H_{1}(C(\bW_{\wdt S_0},\Z))$$ 
is obviously one-free. 

 \bigskip

 We are interested in a slightly more general situation. For any two given subset $\wdt S_h, S^k$ such that $\sharp (\wdt S_h \cup S^k) \leq m$, we consider the complexes 
$F_m^{\wdt S_h \cup S^k}(\wdt \bW)$ generated by cells $E([w],\Gamma)$ such that 
$\Gamma \supset \wdt S_h \cup S^k$.
Hence we define the inclusion maps:
 \begin{equation*}
  i_{m}^h[l]:=i: \bigoplus_{j=1}^{\wdt m_{k}} F^{h+1}_{k+1}(\bW_{\wdt S_{k}})[l] \lgr F_m^{\wdt S_h \cup S^{k+1}}(\wdt \bW)
 \end{equation*}
 as
 \begin{equation*}
\begin{split}
&i(j.E(w_{\wdt S_{k}}, \wdt S_h \cup \Gamma))=i(W^{\wdt S_{k}} (j).E(w_{\wdt S_{k}}, \wdt S_h \cup \Gamma))=\\
& i(w^{\wdt S_{k}}.E(w_{\wdt S_{k}}, \wdt S_h \cup \Gamma))=
[w^{\wdt S_{k}}].E([w_{\wdt S_{k}}], \wdt S_h \cup \Gamma \cup S^{k+1})=E([w], \wdt S_h \cup \Gamma \cup S^{k+1})
\end{split} 
  \end{equation*}
where $W^{\wdt S_{k}}$ is the subset of $W$ isomorphic 
to $(\wdt W/\Lambda)^{\wdt S_{k}}$, $\wdt m_{k}$ its cardinality  and 
$w^{\wdt S_{k}}=W^{\wdt S_{k}}(j)$ its $j$-th element in a fixed order. 

\bigskip
  
They provide short exact sequences as in (\ref{eqn:suc.es.}) and (\ref{eqn:suc.es.tor}):
\begin{equation}\label{sequno}
0 \lgr \bigoplus_{j=1}^{\wdt m_{k}} F^{h+1}_{k+1}(\bW_{\wdt S_{k}})[l] \lgr F_m^{\wdt S_h \cup S^{k+1}}(\wdt \bW) \lgr F_m^{\wdt S_h \cup S^{k}}(\wdt \bW) \lgr 0 
\end{equation}
\bigskip

If $\sharp (\wdt S_h \cup S^k) = m-1$ then $k=h+1$ and, for $l_h=m-h-1$, we get the last short exact sequence:
\begin{equation*} 
0 \lgr \bigoplus_{j=1}^{\wdt m_{h+1}} F^{h+1}_{h+2} (W_{\wdt S_{h+1}})[l_h-1] \stackrel{i}\lgr
F^{\wdt S_{h} \cup S^{h+2}}_m(\wdt \bW) \stackrel{j}\lgr F^{\wdt S_h \cup S^{h+1}}_{m}(\wdt \bW) \lgr 0.
\end{equation*}
It is a simple remark that
$$H_m(F^{\wdt S_h \cup S^{h+1}}_{m}(\wdt \bW),\Z) = F^{\wdt S_h \cup S^{h+1}}_{m}(\wdt \bW) \simeq \bigoplus_{j=1}^{\wdt m_{h}} F^{h+1}_{h+1}(W_{\wdt S_{h}})[l_h] =\bigoplus_{j=1}^{\wdt m_{h}} H^{h+1}(F^{h+1}_{h+1}(W_{\wdt S_h}),\Z)$$ 
are the free modules generated by $E([w],\wdt S\setminus \{s_{h+1}\}) = 
E([w],\wdt S_h \cup S^{h+1}) \simeq E(w,\wdt S_h)=
w^{\wdt S_{h}}.E(w_{\wdt S_h}, \wdt S_h)$. 

Moreover the map $$\wdt \Delta:  F^{\wdt S_h \cup S^{h+1}}_{m}(\wdt \bW)  \lgr \bigoplus_{j=1}^{\wdt m_{h+1}} F^{h+1}_{h+2} (W_{\wdt S_{h+1}})$$
splits as follows:
\begin{equation*}
\begin{array}{ccc}
 \bigoplus_{j=1}^{\wdt m_h} F^{h+1}_{h+1}(W_{\wdt S_h})[l_h] &\stackrel{\Delta}\lgr 
 &\bigoplus_{j=1}^{\sharp W} C(W_{\emptyset}) \\
\wr \mid&  & \downarrow  \scriptstyle{i} \\
F^{\wdt S_h \cup S^{h+1}}_{m}(\wdt \bW)& \xrightarrow{\wdt \Delta} & \bigoplus_{j=1}^{\wdt m_{h+1}} F^{h+1}_{h+2}(W_{\wdt S_{h+1}}) 
\end{array}
\end{equation*}

and we get the commutative diagram in homology:
\begin{equation*} \begin{array}{ccc}
\bigoplus_{j=1}^{\wdt m_h} H_{h+1}(F^{h+1}_{h+1}(W_{\wdt S_h}), \Z) &\xrightarrow{\Delta_*} &  \bigoplus_{j=1}^{\sharp W} H_{0}(C(\bW_{\emptyset}),\Z) \\
\wr \mid & &\downarrow\scriptstyle{i_*} \\
H_{m}(F^{\wdt S_h \cup S^{h+1}}_{m}(\wdt \bW),\Z)& \xrightarrow{\wdt \Delta_*} & \bigoplus_{j=1}^{\wdt m_{h+1}}H_{h+1} (F^{h+1}_{h+2}(W_{\wdt S_{h+1}}) ,\Z) .
\end{array}
\end{equation*}

Hence if $\sharp(\wdt S_h \cup S^k)=m-1$ the map $\wdt \Delta_*$ is one-free since it is composition of one-free maps $\Delta_*$ and $i_*$. So far we proved the base of a more general induction.

Going backwards on homology exact sequences induced by (\ref{sequno}) we get maps
\begin{equation}\label{deltastar}
\wdt \Delta_*: H_{*+1}(F^{\wdt S_h \cup S^k}_{m}(\wdt \bW), \Z) \lgr \bigoplus 
H_{*-l}(F^{h+1}_{k+1}(W_{\wdt S_k}), \Z).  
\end{equation}
Let us assume, by induction, that they are one-free maps for all $\wdt S_h, S^k$ such that $n < \sharp (\wdt S_h \cup S^k) \leq m-1$ (i.e. $H_*(F^{\wdt S_h \cup S^{k}}_{m}(\wdt \bW), \Z)$ are free modules for $n \leq \sharp (\wdt S_h \cup S^k) \leq m-1$). 

Let  $\sharp (\wdt S_h \cup S^k)$ be equal to  $n$.

\bigskip

We can also filter $F^{\wdt S_h \cup S^k}_{m}(\wdt \bW)$ as follows:
\begin{equation*}
\begin{split}
&i^m [h+1]:=i: \bigoplus_{j=1}^{m^1\ldots m^{h+1}} F_{l_h}^{m-k}(\bW_{S^{h+1}})[h+1] \lgr  F^{\wdt S_h \cup S^k}_{m}(\wdt \bW)\\
&i(w^{S^{h+1}}E(w_{S^{h+1}}, \Gamma \cup S^{k}))= E([w], \wdt S_h \cup \Gamma \cup S^k).
\end{split}
\end{equation*}

\bigskip

We get the exact sequences
\begin{equation*}
0 \lgr \bigoplus_{j=1}^{m^1\ldots m^{h+1}} F_{l_h}^{m-k}(\bW_{S^{h+1}})[h+1] \lgr  F^{\wdt S_h \cup S^k}_{m}(\wdt \bW) \lgr
F^{\wdt S_{h+1} \cup S^k}_{m}(\wdt \bW) \lgr 0. 
\end{equation*}

This is equivalent to say that for any cell $E([w],\wdt S_h \cup \Gamma \cup S^k) \in F^{\wdt S_h \cup S^k}_{m}(\wdt \bW)$ we have only two possibilities: 
\begin{equation*}
\begin{split}
&i) s_{h+1} \in \Gamma \mbox{ and hence } E([w],\wdt S_h \cup \Gamma \cup S^k) =E([w],\wdt S_{h+1} \cup \Gamma^{\prime} \cup S^k) \in  F^{\wdt S_{h+1} \cup S^k}_{m}(\wdt \bW) \\ 
&\mbox{or} \\ 
&ii) s_{h+1}\notin \Gamma \mbox{ and hence }
E([w],\wdt S_h \cup \Gamma \cup S^k)=i(w^{S^{h+1}}E(w_{S^{h+1}}, \Gamma \cup S^k)) \in \\
&\in i(\bigoplus_{j=1}^{m^1\ldots m^{h+1}} F_{l_h}^{m-k}(\bW_{S^{h+1}})[h+1]).
\end{split}
\end{equation*}
As a consequence if $\wdt \Delta:F^{\wdt S_h \cup S^k}_{m}(\wdt \bW) \lgr 
\bigoplus_{j=1}^{\wdt m_k} F^{h+1}_{k+1}(W_{\wdt S_k})[l] $ is the map which induces the map 
$\wdt \Delta_*$ in (\ref{deltastar}),  
$\wdt \Delta$ splits as follows:  
\begin{displaymath}
\left[
\begin{array}{cc}
\wdt \Delta_{\mid_{F_m^{\wdt S_{h+1} \cup S^k}(\wdt W)}} & 0 \\
0  & \Delta_{\mid_{\bigoplus F_{l_h}^{m-k}(\bW_{S^{h+1}})}}
\end{array}
\right]_.
\end{displaymath}
Here $\wdt \Delta_{\mid_{F_m^{\wdt S_{h+1} \cup S^k}(\wdt W)}}$ is the map $\wdt \Delta$
defined on  $F_m^{\wdt S_{h+1} \cup S^k}(\wdt W)$, i.e. on a complex such that 
$\sharp(\wdt S_{h+1} \cup S^k)=n+1$ if $\sharp(\wdt S_{h} \cup S^k)=n$. 

From now on we will denote this map $\wdt \Delta_{n+1}$ in order to distinguish it from 
$\wdt \Delta_{n}$.

\bigskip

By previous consideration it follows that the diagram on complexes
\begin{equation}\label{doppieq} 
\begin{array}{cccccc}
0 \lgr \bigoplus F_{l_h}^{m-k}(\bW_{S^{h+1}})[h+1] &\xrightarrow{\wdt i} &  F^{\wdt S_h \cup S^k}_{m}(\wdt \bW) &\xrightarrow{\wdt j} & F^{\wdt S_{h+1} \cup S^k}_{m}(\wdt \bW)  & \lgr 0 \\
\Delta \downarrow & & \wdt \Delta_n \downarrow &  &\wdt \Delta_{n+1} \downarrow & \\
0 \lgr \bigoplus \bigoplus C(\bW_{\wdt S_{k} \setminus \wdt S_{h+1}})[l][h+1] &\xrightarrow{i} & 
\bigoplus F^{h+1}_{k+1}(W_{\wdt S_k})[l]  &\xrightarrow{j} &  
\bigoplus F^{h+2}_{k+1} (W_{\wdt S_k})[l]  & \lgr 0
\end{array}
\end{equation}
is commutative.  

Here $i:\bigoplus \bigoplus C(\bW_{\wdt S_{k} \setminus \wdt S_{h+1}})[l][h+1] \lgr 
\bigoplus_{j=1}^{\wdt m_k} F^{h+1}_{k+1}(W_{\wdt S_k})[l] $ is the map of type (\ref{eqn4}) such that 
$i(w^{\wdt S_{k} \setminus \wdt S_{h+1}}.E(w_{\wdt S_{k} \setminus \wdt S_{h+1}},\Gamma))=
w^{\wdt S_k}E(w_{\wdt S_k},\wdt S_h \cup \Gamma)$.

\bigskip

Let us remark that the sum 
$$\bigoplus \bigoplus C(\bW_{\wdt S_{k} \setminus \wdt S_{h+1}})[l][h+1]=\bigoplus_{j=1}^{\sharp W /\sharp W_{\wdt S_{k} \setminus \wdt S_{h+1}}}C(\bW_{\wdt S_{k} \setminus \wdt S_{h+1}})[l][h+1]$$
splits in different ways depending if we are considering the horizontal exact sequence or the 
vertical map $\Delta$. 

The diagram (\ref{doppieq}) gives rise to the following commutative diagram in homology:

\begin{equation*} 
\begin{array}{cccccc}
\lgr \bigoplus H_{*-h-1}(F_{l_h}^{m-k}(\bW_{S^{h+1}}),\Z) &\xrightarrow{\wdt i_*} &  H_*(F^{\wdt S_h \cup S^k}_{m}(\wdt \bW), \Z) &\xrightarrow{\wdt j_*} & H_*(F^{\wdt S_{h+1} \cup S^k}_{m}(\wdt \bW), \Z)  & \lgr \\ 
\Delta_* \downarrow & & \wdt \Delta_{n *} \downarrow &  &\wdt \Delta_{n+1 *} \downarrow & \\
\lgr \bigoplus H_{*-l-h-1}(C(\bW_{\wdt S_k \setminus \wdt S_{h+1}}),\Z) &\xrightarrow{i_*} & 
\bigoplus H_{*-l}(F^{h+1}_{k+1}(W_{\wdt S_k}), \Z)  &\xrightarrow{j_*} &  \bigoplus H_{*-l}(F^{h+2}_{k+1}(W_{\wdt S_k}), \Z)   & \lgr
\end{array}
\end{equation*}

\bigskip

The maps $i_*$, $j_*$ and $\Delta_*$ are one-free (see section \ref{sezione2}).

Moreover $H_{*-h-1}(F_{l_h}^{m-k}(\bW_{S^{h+1}}),\Z)$ and
$H_*(F^{\wdt S_h \cup S^k}_{m}(\wdt \bW), \Z)$ are free modules respectively by theorem \ref{teofmk} and by inductive hypothesis. Then the maps  $\wdt i_*$ and $\wdt j_*$ in the diagram are one-free. Moreover $\wdt \Delta_{n+1 *}$ are one-free by induction and hence we get that maps $\wdt \Delta_{n*}$ are one-free too.

So far we proved the main result of the paper:

\begin{teo}The integer (co)-homology of the complement $\Cal R_W$  is torsion free.
\end{teo}

As an immediate consequence of the above theorem, $H^*(\Cal R_W, \Z)$ coincides with the De Rham cohomology described in \cite{DP} and the Betti numbers can be easily computed 
using results in \cite{Mo1}.

\bigskip

In general we have the following

\begin{conj}
Let  $\mathcal{T}_{X}$ be a thick toric arrangement in the sense of \cite{MoSe}. Then the integer cohomology of the complement is torsion free (and hence it coincides with the De Rham cohomology computed in \cite{DP}).
\end{conj}

\end{document}